\documentclass [12pt,a4paper,reqno]{amsart}
\usepackage{amsmath,amsthm,amscd}
\usepackage{amsfonts,amssymb}
\allowdisplaybreaks

\DeclareMathOperator{\esssup}{esssup}

\newcommand{\be}{\begin{equation}}
\newcommand{\ef}{\end{equation}}
\chardef\bslash=`\\ 

\newcommand{\wt}{\widetilde}
\newcommand{\wh}{\widehat}
 \renewcommand{\sectionmark}[1]{}
\renewcommand{\Im}{\operatorname{Im}}

\newcommand{\ve}{\varepsilon}

\newcommand{\bk}{\bigskip}
\newcommand{\iy}{\infty}

\newcommand{\const}{\operatorname{const}}

 \usepackage{amsfonts}
\newcommand{\field}[1]{\mathbb{#1}}

\newcommand{\dl}{\delta}
\newcommand{\D}{\field D}
\newcommand{\om}{\omega}
\newcommand{\z}{\zeta}
\newcommand{\ov}{\overline}
\newcommand{\vp}{\varphi}
\newcommand{\hC}{\widehat{\field{C}}}
\newcommand{\C}{\field{C}}
\newcommand{\R}{\field{R}}
\newcommand{\B}{\mathbf{B}}
\newcommand{\T}{\mathbf{T}}
\newcommand{\fH}{\field{H}}

\newcommand{\fc}{\frac}

\newcommand{\dist}{\operatorname{dist}}

\newcommand{\Belt} {\operatorname{Belt}}

\newcommand{\vk} {\varkappa}
\newcommand{\x} {\mathbf x}
\renewcommand{\a} {\alpha}

\newcommand{\ld}{\lambda}

\begin{document}

\title{Extremal quasiconformality vs rational approximation}

\author{Samuel L. Krushkal}

\begin{abstract} We show that on most of the hyperbolic simply
connected domains the weighted bounded rational approximation in a
natural sup norm is possible only for a very sparse set of
holomorphic functions (in contrast to integral approximation). The
obstructions are caused by the features of extremal
quasiconformality.
\end{abstract}


\maketitle

\bigskip

{\small {\textbf {2010 Mathematics Subject Classification:} Primary:
30C62, 30C75, 30E10; Secondary: 30F45, 30F60, 32G15}

\medskip

}

\bigskip

\markboth{S. L. Krushkal}{Extremal quasiconformality vs bounded
rational approximation} \pagestyle{headings}

\bigskip
The paper is dedicated to the 100th anniversary of  
Georgii Dmitrievich Suvorov, my first university adviser and teacher. 
He was an outstanding mathematician and a widely talented,  
extremely great human being.

\bigskip\bigskip
\centerline{\bf 1. RESULTS}

\bigskip

This paper gives a link of geometric function theory to weighted
bounded rational approximation of holomorphic functions in sup norms and shows how the intrinsic features of extremal quasiconformal maps
and universal Teichm\"{u}ller space provide strong obstructions to such an approximation. The situation is completely different from the
integral approximation.

\bigskip\noindent {\bf 1.1. Introductory remarks}. The classical
directions in rational interpolation of holomorphic functions
investigated by many authors concern mainly the uniform
interpolation of functions holomorphic in the inner points of the closed
set $X$ on the Riemann sphere $\hC = \C \cup \{\iy\}$ and continuous
on $X$ by rational functions with poles off $X$ and its combination with interpolation (see, e.g., \cite{Za}). The second approach was originated by Walsh (see
\cite{Wa}, \cite{AmW}) and has  recently been extended in \cite{Gu} to the
functional space $A^{-\iy} = \bigcup_{q>0} A^{-q}$ over the Dini
domains (with topology of the inductive limit), where $A^{-q}(D)$ is the Banach space of holomorphic functions in a domain $D$ with norm $\|f\|= \sup_D \dl_D(z) |f(z)|$; 
here $\dl_D(z) = \dist(z, \partial D)$ denotes the Euclidean distance
from the point $z \in D$ to the boundary.

More generally, one considers a space $\mathcal F$ of holomorphic \
functions $f$ in a domain $D \subset \hC$. For any given collections
of points
$$
A_n = \{a_{nj}\}_{j=0}^n \subset D, \quad B_n = \{b_{nj}\}_{j=1}^n
\subset \hC \setminus D \quad (n = 1, 2, \dots),
$$
there exists a unique rational function $r_{n,f}$ of degree $n$,
with poles at $B_n$, interpolating to $f$ at $A_n$, counting
multiplicities. The problem is to select these collections $A_n$ and
$B_n$ so that for all $f$ the interpolants $r_{n, f}$ converge to
$f$ in the topology of $\mathcal F$.

It is established in \cite{Gu} (by sweeping out appropriate
measures and application of the potential methods) that for $f \in A^{-q}(D), \ q > 0$, the interpolants $r_{n, f}$ are convergent
(under appropriate conditions) to $f$ in $A^{-q^\prime}$ norm, where $q^\prime(q) \gg q$.

\bigskip\noindent {\bf 1.2. Weakened rational approximation in sup
norm}. Let $L$ be an oriented quasiconformal Jordan curve
(quasicircle) on the Riemann sphere $\hC = \C \cup \{\iy\}$ with the
interior and exterior domains $D$ and $D^*$, and let $p \ge 2$ be an
integer. Denote by $\ld_D(z) |dz|$ the hyperbolic metric of $D$ of
Gaussian curvature $- 4$ and consider the Banach spaces $A_p(D), \
\B_p(D)$ of holomorphic functions (quadratic differentials) $\vp$
with norms
$$
\|\vp\|_{\B_p} = \sup_D \ld_D(z)^{-p} |\vp(z)|, \quad \|\vp\|_{A_p}
= \iint\limits_D \ld_D(z)^{2-p} |\vp(z)| dx dy,
$$
respectively; due to \cite{Be2}, the space $\B_p$ is dual to $A_p$.

Note that for simply connected domains (with more than one boundary
points) not containing inside the infinite point,
 \be\label{1}
\fc{1}{4} \le \ld_D(z) \dl_D(z) \le 1,
\end{equation}
where the right hand inequality follows from the Schwarz lemma and
the left from Koebe's $\fc{1}{4}$ theorem (so the spaces $A^{-p}$
mentioned above are obtained by renormalization of $\B_p$).

We also shall use the notations
$$
\D = \{z: \ |z| < 1\}, \quad \D^* = \{z \in \hC: \ |z| > 1\}; \quad
\fH = \{z: \Im z > 0\}, \ \fH^* = \{z: \Im z < 0\}.
$$

\bigskip
We start with the following general

\bigskip\noindent
{\bf Theorem 1}. {\em Let $D \subset \hC$ be a domain with
quasiconformal boundary $L$. Then for any function $f \in \B_p(D)$
there exists a sequence of rational functions with poles of order
two on $L$ of the form
 \be\label{2}
r_n(z) = \sum\limits_1^n \fc{c_j}{(z - a_j)^2}, \quad
\sum\limits_1^n |c_j| > 0,
\end{equation}
such that $\lim\limits_{n\to \iy} \|r_n - \vp\|_{\B_{p+1}(D)} = 0$.}

\bigskip In the case of the half-plane (or disk), this theorem is
strengthened as follows.

\bigskip\noindent
{\bf Theorem 2}. {\em For any $\vp \in \B_p(\fH)$ there exists a
sequence of rational functions (2) with real poles $a_j$ and real
coefficients $c_j$, convergent to $\vp$ in $\B_{p+1}(\fH)$. }

\bigskip
One can see from the proof of Theorem 1 that the weight exponent
$p + 1$ is not sharp; though it is not clear
whether this exponent can be replaces by $p + \epsilon(p)$ with $0 <
\epsilon(p) < 1$.

The limit case $\epsilon = 0$ has intrinsic interest. Then the
assertion on convergence fails for $\epsilon = 0$, because the space
$\B_2$ of bounded holomorphic quadratic differentials $\vp dz^2$ is
not separable for any Riemann surface of infinite genus, and similar
for all $\B_p$.

This was established by functional-analytic methods but can be
established also from geometric features generated by  Thurston's
theorem on existence of uncountable many conformally rigid domains
(see \cite{Th}, \cite{As}). Such domains correspond to the isolated
points of the set $\mathbf U \setminus \T$ in $\B_2(\fH)$, where
$\mathbf U$ denotes the set of the Schwarzian derivatives
$$
S_w(z) = \Bigl(\fc{w^{\prime\prime}(z)}{w^\prime(z)}\Bigr)^\prime -
\fc{1}{2} \Bigl(\fc{w^{\prime\prime}(z)}{w^\prime(z)}\Bigr)^2
$$
of all univalent functions on $\fH$ and $\T$ is the universal
Teichm\"{u}ller space modeled by a bounded domain in $\B_2(\fH)$
(formed by $f$ having quasiconformal extension to $\hC$).

Note also that any $\vp \in \B_2(D)$ can be regarded as the
Schwarzian of a locally univalent function in $D$ and determines
this function up to a Moebius transformation of $\hC$.

\bigskip\noindent{\bf 1.3. Main theorems}. Our aim is to show that on most of the hyperbolic simply connected domains the rational approximation in $\B_2$  is possible only for a very sparse subset of functions.

\bigskip\noindent
{\bf Theorem 3}. {\em For any simply connected domain $D \subset
\hC$ with quasiconformal boundary $L$, whose conformal mapping function $g_D: \fH \to D$ satisfies $\|S_{g_D}\|_{\B_2} < 1/2$, the set of functions $\vp \in
\B_2(D)$ approximated in $\B_2$ norm by general rational functions
with poles of order two on $L$,
  \be\label{3}
r_n(z) = \sum\limits_1^n \fc{c_j}{(z - a_j)^2} + 
\sum\limits_1^n \fc{c_j^\prime}{z - a_j}, \quad
\sum\limits_1^n |c_j| > 0,
\end{equation}
is nonwhere dense in the space $\B_2(D)$. }

\bigskip
This theorem is a consequence of some deep results concerning the extremal  maps and the Grunsky operator given by

\bigskip\noindent
{\bf Theorem 4}. {\em There exists a constant $c_0 > 0$ 
such that for any simply connected domain $D \subset \hC$ with quasiconformal boundary $L$ and such that $\|S_{g_D}\|_{\B_2} < 1/2$, and for any rational function $r_n$ with poles of order two on $L$ of the form (3) and norm
$\|r_n\|_{\B_2(D)} < c_0$, we have the equalities 
 \be\label{4}
\vk_D(w) = k(w) = \|r_n\|_{\B_2(D)},
\end{equation}
where $\vk_D(w)$ and $k(w)$ denote the Grunsky and Teichm\"{u}ller
norms of (appropriately normalized) univalent solution $w: \ D \to \hC$ of the Schwarzian equation $S_f = r_n$.} 

Note that the indicated constant $c_0$ does not depend on $D$. In the case of a disk (half-plane), one can take $c_0 = 1/2$. 

\bigskip
Theorem 4 shows that all Schwarzians $\vp = r_n$ with
$\|r_n\|_{\B_2} < c_0$ are not Strebel points in the universal
Teichm\"{u}ller space $\T$ (in other words, the conformal maps $f$
with $S_f = r_n$ do not have Teichm\"{u}ller extremal extensions onto the complementary domain $D^* = \hC \setminus \ov D$). 
Hence, such points cannot be dense in $\B_2(D)$. 

Note also that, in view of the
first equality in (4), the dilatation $k(f)$ is attained on the
squares of holomorphic abelian differentials $\om dz$ on $D$.

\bigskip
Theorem 3 is an immediate consequence of the equalities (4) in view of either of two basic results on openness and density in the
universal Teichm\"{u}ller space $\T$: the result of \cite{Kr3} for
points with non-equal Teichm\"{u}ller and Grunsky norms and Lakic's result \cite{La} on the density of Strebel points (in arbitrary Teichm\"{u}ller space). 

Indeed, it suffices to establish the assertion of Theorem 3 for $\vp \in \B_2(D)$ with $\|\vp\| < c_0$. 
Both Teichm\"{u}ller and Grunsky norms are continuous on the universal Teichm\"{u}ller space $\T$ modelled as a bounded domain in $\B_2$ containing the origin. Hence, the equality (4) is preserving also for the limit function 
$f = \lim r_n$ of any sequence of rational functions in $\B_2$. 

But this equality implies that any $f \in \B_2$ 
with sufficiently small $\|f\|$ generates the \linebreak Beltrami coefficient 
$$
\mu(z) = \fc{1}{2} |z - \ov z|^2 f \circ g_D (\ov z), 
$$
where $g_D$ is a conformal map of $\fH$ onto domain $D$ (called harmonic), and this coefficient is extremal 
in its equivalence class. On the other side, it is not of Teichm\"{u}ller type. 

The latter is impossible in view of the indicated above openness results.

\bigskip\noindent{\bf 1.4. Generalization of Theorems 3 and 4}. One
can see from the proof of Theorem 4 in Section 5 that actually its
arguments are valid for an arbitrary meromorphic function $\vp$ on $\C$
having poles of order two which are located on a quasicircle $L$
passing through the infinite point and accumulate to this point.
This gives the following extension of the above theorems.

\bk\noindent
{\bf Theorem 5}. {\it For any simply connected domain $D \in \hC$ with quasiconformal boundary $L$ passing through the infinite point whose conformal mapping function satisfies $\|S_{g_D}\|_{\B_2} < 1/2$, the subspace $\mathcal M_2$ in $\B_2(D)$ formed by
meromorphic functions $\vp$ on $\C$ with poles of order two, which are located on $L$ and accumulate to $\iy$, is nonwhere dense in the space $\B_2(D)$.

If $\|\vp\|_{\B_2(D)}$ is sufficiently small, so that the Schwarzian equation $S_w = \vp$ has a univalent solution $w(z)$ on $D$, then}
$$
\vk_D(w) = k(w) = \|\vp\|_{\B_2(D)}.
$$

\bigskip\bigskip
\centerline{\bf 2. PROOF OF THEOREM 1}

\bigskip

First note that $A_p(D) \subset \B_p D)$, and for any $\vp \in
A_p(D)$,
  \be\label{5}
\|\vp\|_{\B_p} \le \fc{4}{\pi} \|\vp\|_{A_p}.
\end{equation}
Indeed, since both these norms are conformally invariant, it
suffices to verify (5) for $D = \fH$. Then $\ld_D (z) = 1/(2 y)$,
and applying the mean inequality for holomorphic \ functions, one
gets
$$
|f(z)| \le \fc{1}{\pi y^2} \iint\limits_{|z- \z| \le y} |f(\z)| d\xi
d\eta \le \fc{(2y)^{2-p}}{\pi y^2} \iint\limits_{|z- \z| \le y} (2
\eta)^{2-p} d\xi d \eta
$$
(here $\z = \xi + i \eta, \ z \in D, \eta \le 2 y$), which yields
(5). It follows also that for any $\vp \in A_p(D)$,
 \be\label{6}
\lim\limits_{\rho \to 0} \sup_{\dl_D(z)\ge \rho} \ld_D^p(z) |\vp(z)|
= 0.
\end{equation}

Without loss of generality, one can assume that the boundary curve
$L \ni \iy$ and $0 \in D$. For any $\vp \in \B_p(D)$,
$$
|\vp(z)| \le \|\vp\|_{B_p} \ld_D(z)^{-p} \asymp \|\vp\|_{B_p}
\dl_D(z)^{-p};
$$
hence it belongs to $A_{p+1}$ and its integral
$$
I_\vp(z) = \int_0^z
\vp(\z) d \z
$$
belongs to $A_p$.

By the Bers approximation theorem \cite{Be1}, there exists a
sequence of rational functions $\wt r_j(z)$ with simple poles on $L$
and no other singularities, such that\footnote{This theorem is proved in \cite{Be1} for the integrable holomorphic functions $f \in A_2$; the proof for the weighted spaces $A_p$ can be done along the same lines, see \cite{Kr1}.}
 \be\label{7}
\lim\limits_{n \to \iy} \|\wt r_j - I_\vp\|_{A_p} = 0.
\end{equation}
Further, due to \cite{Be2}, for every $\psi \in \B_p(D)$ the
following reproducing formula is valid: 
  \be\label{8}
\psi(z) = - \fc{2p - 1}{\pi} \iint\limits_D \fc{(\z -
h(\z)^{2p-2}(\partial h(\z)/\partial \ov \z)) \psi(h(\z))}{(\z -
z)^{2p}} d\xi d \eta;
\end{equation}
here $\z \mapsto h(\z)$ is a quasiconformal reflection with respect to the quasicircle $L = \partial D$ (i.e., orientation reversing
quasiconformal automorphism of $\hC$ mapping $D$ onto its
complementary domain and leaving fixed all points of $L$) which is
uniformly bilipschizian on $\C$, i.e., for all points $z_1, z_2 \in
\C$ the inequality
$$
c_0^{-1} |z_1 - z_2| \le |h(z_1) - h(z_2)| \le c_0 |z_1 - z_2|
$$
holds with some constant $c_0 > 1$. Moreover, the numerator
$$
\nu_\psi(\z) = - \fc{2p - 1}{\pi} (\z - h(\z))^{2p-2} \fc{\partial
h(\z)}{\partial \ov \z}
$$
is estimated uniformly by
$$
|\nu_\psi(\z)| \le c(c_0) \|\psi\|_{\B_2(D)} \ld_D(\z)^{2-p}.
$$
Applying (8) to $\psi = \wt r_j - I_\vp$ and differentiating both
sides in $z$, one obtains from (7) that $r_j = \wt r_j^\prime$ are
convergent to $\vp$ in $\B_{p+1}(D)$, completing the proof of the theorem.

\bigskip\bigskip
\centerline{\bf 3. PROOF OF THEOREM 2}

\bigskip
We apply the following result on integral approximation in the unit
disk given in \cite{Kr1} improving for the disk the Bers
approximation theorem mentioned above and also related to the theory
of extremal quasiconformal maps.

\bigskip\noindent{\bf Proposition 1}. {\em Let $p$ and $m$ be two
integers such that $p \ge 2$ and $m \ge 1$. Then for any $\psi \in
A_p(\D)$, there exists a sequence of rational functions $\wt r_j$
which have only simple poles on the unit circle $S^1$ and satisfy
the condition $\Im [\z^m \wt r_j(\z)] = 0$ on $S^1$ (outside of the poles of $\wt r_j$), such that}
$$
\lim\limits_{j \to \iy} \|\wt r_j - \psi\|_{A_p} = 0.
$$

For even $m$ and $\z = e^{i \theta} \in S^1$, we have
 \be\label{9}
\Im [\wt r_j(\z) d \z^m] = (- 1)^{m/2} \Im [\z^m \wt r_j(\z)] d
\theta^m;
\end{equation}
thus the above proposition can be reformulated as follows.

\bigskip\noindent{\bf Proposition 2}. {\em  For any function
$\psi \in A_p(\D)$, there exists a sequence of rational functions
$\wt r_j$ which have only simple poles on the circle $S^1$ and
satisfy the condition $\Im [\wt r_j(\z) d \z^{2m}] = 0$ on $S^1$,
such that $\lim\limits_{j \to \iy} \|\wt r_j - \psi \|_{A_p} = 0$.}

\bigskip
Take $m = p$ and, applying the fractional linear map $\sigma(z) = (z - i)/(z + i)$ of the upper half-plane onto the unit disk,
approximate similar to Theorem 1 the integrated functions
$$
\sigma_{*} I_\psi = (I_\psi \circ \sigma) (\sigma^\prime)^{4m-2}
$$
by the corresponding rational functions
$$
\sigma_{*} \wt r_j = (\wt r_j \ \circ \sigma) (\sigma^\prime)^{4m-2} \in
A_{p+1}(\fH)
$$
having real poles and coefficients in view of (9). Now applying the
reproducing formula for the upper half-plane,
$$
\psi(z) =  \fc{2p - 1}{\pi} \iint\limits_D \fc{(\z - \ov
\z)^{2p-2}\psi(\ov \z)}{(\z - z)^{2p}} d\xi d \eta,
$$
one straightforwardly obtains the conclusion of Theorem 2.

\bigskip\bigskip
\centerline{\bf 4. BACKGROUNDS OF THEOREM 3}

\bigskip

As was mentioned, Theorem 3 relies on some intrinsic features of extremal quasiconformal maps and the Grunsky operator. For convenience, we briefly describe here these underlying results. 

\bigskip\noindent
{\bf 4.1. Extremal quasiconformality}. 
Let $L$ be a quasicircle passing through the points $0, 1, \iy$
which is the common boundary of two domains $D$ and $D^*$. Take the unit ball of Beltrami coefficients supported on $D^*$,
$$
\Belt(D^*)_1 = \{\mu \in L_\iy(\C): \ \mu|D = 0\, \ \ \|\mu\|_\iy < 1\}
$$
and consider the corresponding quasiconformal automorphisms
$w^\mu(z)$ of the sphere $\hC$ satisfying on $\C$ the  Beltrami equation $\ov{\partial} w = \mu \partial w$ preserving the points $0, 1, \iy$ fixed. We call the quantity $k(w) = \|\mu_w\|_\iy$ the {\bf dilatation} of the map $w$.

Take the equivalence classes $[\mu]$ and $[w^\mu]$ letting the coefficients $\mu_1$ and $\mu_2$ from $\Belt(D^*)_1$ be equivalent if
the corresponding maps $w^{\mu_1}$ and $w^{\mu_2}$ coincide on $L$ (and hence on $\ov{D}$). These classes are in one-to-one correspondence with the Schwarzians $S_{w^\mu}$ on $D$ which fill a
bounded domain in the space $\B_2(D)$ modelling the universal
Teichm\"{u}ller space $\T = \T(D)$ with the base point $D$. The quotient map
$$
\phi_\T: \ \Belt(D^*)_1 \to \T, \quad \phi_\T(\mu) = S_{w^\mu}
$$
is holomorphic (as the map from $L_\iy(D^*)$ to $\B_2(D)$). Its intrinsic {\it Teichm\"{u}ller metric} is defined by
$$
\tau_\T (\phi_\T (\mu), \phi_\T (\nu)) = \frac{1}{2} \inf \bigl\{
\log K \bigl( w^{\mu_*} \circ \bigl(w^{\nu_*} \bigr)^{-1} \bigr) : \
\mu_* \in \phi_\T(\mu), \nu_* \in \phi_\T(\nu) \bigr\},
$$
It is the integral form of the infinitesimal Finsler metric
$$
  F_\T(\phi_\T(\mu), \phi_\T^\prime(\mu) \nu) = \inf
\{\|\nu_*/(1 - |\mu|^2)\|_\iy: \ \phi_\T^\prime(\mu) \nu_* =
\phi_\T^\prime(\mu) \nu\}
$$
on the tangent bundle $\mathcal T\T$ of $\T$, which is locally
Lipschitzian.

We call the Beltrami coefficient $\mu \in \Belt(D^*)_1$ {\bf extremal} (in its class) if
$$
\|\mu\|_\iy = \inf \{\|\nu\|_\iy: \ \ \phi_\T(\nu) = \phi_\T(\mu)\}
$$
and call $\mu$ {\bf infinitesimally extremal} if
$$
\|\mu\|_\iy = \inf \{\|\nu\|_\iy: \ \ \nu \in L_\iy(D^*), \ \
\phi_\T^\prime(\mathbf 0) \nu = \phi_\T^\prime(\mathbf 0) \mu\}.
$$
Any infinitesimally extremal Beltrami coefficient $\mu$ is globally extremal (and vice versa), and by the basic
Hamilton-Krushkal-Reich-Strebel theorem the extremality of $\mu$ is equivalent to the equality
$$
\|\mu\|_\iy = \inf \{|<\mu, \psi>_{D^*}|: \ \ \psi \in A_2(D^*): \
\|\psi\| = 1\}
$$
(where $A_2(D^*)$ is the subspace of $L_1(D^*)$ formed by
holomorphic \ functions on $D^*$) and the pairing
$$
\langle \mu, \psi \rangle_{D^*} = \iint_{D^*}  \mu(z) \psi(z) dx dy, \quad \mu \in L_\iy(D^*), \ \psi \in L_1(D^*) \ \ (z = x + iy).
$$

Let $w_0 :=w^{\mu_0}$ be an extremal representative of its class
$[w_0]$ with dilatation
$$
k(w_0) = \|\mu_0\|_\iy = \inf \{k(w^\mu): w^\mu|L = w_0|L\},
$$
and assume that there exists in this class a quasiconformal map
$w_1$ whose Beltrami coefficient $\mu_{A_1}$ satisfies the
inequality $\esssup_{A_r} |\mu_{w_1}(z)| < k(w_0)$ in some ring
domain $\mathcal R = D^* \setminus G$ complement to a domain $G
\supset D^*$. Any such $w_1$ is called the {\bf frame map} for the
class $[w_0]$, and the corresponding point in the universal
Teichm\"{u}ller space $\T$ is called the {\bf Strebel point}.

\bigskip
These points have the following important properties.

\bigskip\noindent{\bf Proposition 3}. {\em (i) \ If a class $[f]$ has a frame map, then the extremal map $f_0$ in this class (minimizing the dilatation $\|\mu\|_\iy$) is unique and either a conformal or a Teichm\"{u}ller map with Beltrami coefficient $\mu_0
= k |\psi_0|/\psi_0$ on $D^*$, defined by an integrable holomorphic
quadratic differential $\psi_0$ on $D^*$ and a constant $k \in (0,
1)$ \cite{St}.

(ii) \ The set of Strebel points is open and dense in $\T$ 
\cite{La}, \cite{GL}.}

\bigskip

The first assertion holds, for example, for asymptotically conformal
(hence for all smooth) curves $L$. Similar results hold also for
arbitrary Riemann surfaces (cf. \cite{EL}, \cite{GL}).

The boundary dilatation $H(f)$ admits also a local version $H_p(f)$
involving the Beltrami coefficients supported in the neighborhoods
of a boundary point $p \in \partial D$. Moreover (see, e.g., [8, Ch.
17]), $H(f) = \sup_{p\in \partial D} H_p(f)$, and the points with
$H_p(f) = H(f)$ are called {\bf substantial} for $f$ and for its
equivalence class.

\bigskip\noindent
{\bf 4.2. The Grunsky-Milin inequalities}. 
Let $D^* \ni \iy$ be a simply connected domain with quasiconformal
boundary and $\Sigma^0(D^*)$ denote the class of univalent
$\hC$-holomorphic \ functions in $D^*$ with expansions $f(z) = z +
b_0 + b_1 z^{-1} + \dots$ near $z = \iy$ admitting quasiconformal extensions to $\hC$. Their Grunsky-Milin coefficients $\a_{mn}$ are defined from the expansion
 \be\label{10}
- \log \fc{f(z) - f(\z)}{z - \z} = \sum\limits_{m, n = 1}^\iy
\fc{\a_{m n}}{\chi(z)^m \ \chi(\z)^n},
\end{equation}
choosing the branch of the logarithmic function which vanishes as $z = \z \to \iy$. Here $\chi$ denotes a conformal map of $D^*$ onto the disk $\D^*$ so that $\chi(\iy) = \iy, \ \chi^\prime(\iy) > 0$.

Each coefficient $\a_{m n}(f)$ in (10) is a polynomial of a finite number of the initial coefficients $b_1, b_2, \dots, b_{m+n-1}$ of $f$; hence it depends holomorphically on Beltrami coefficients of
extensions of $f$ as well as on the Schwarzian derivatives $S_f \in
\B_2(D^*)$.

A theorem of Milin extending the Grunsky univalence criterion for the disk $\D^*$ states that a holomorphic \ function $f(z) = z + \const + O(z^{-1})$ in a neighborhood of $z = \iy$ can be continued to a univalent function in the whole domain $D^*$ if and only if the
coefficients $\a_{mn}$ satisfy the inequality
$$
\Big\vert\sum\limits_{m,n = 1}^{\iy} \ \sqrt{m n} \ \a_{m n} x_m x_n \Big\vert \le 1
$$
for any point $\x = (x_n)$ from the unit sphere $S(l^2)$ of the Hilbert space of sequences $\mathbf x = (x_n)$ with $\|\mathbf x \|^2 = \sum\limits_{1}^{\iy} |x_n|^2$ (cf. \cite{Gr}, \cite{Mi}, \cite{Po}). We
call the quantity
$$
\vk_{D^*}(f) = \sup \Big\{ \Big\vert \sum\limits_{m,n = 1}^{\iy} \ \sqrt{m n} \ \a_{mn} \ x_m x_n \Big\vert : \ {\mathbf x} = (x_n) \in S(l^2)\Big\} 
$$
the {\bf Grunsky norm} of $f$. The inequality $\vk_{D^*}(f) \le 1$
is necessary and sufficient for univalence of $f$ in $D^*$ (see \cite{Gr}, \cite{Mi}, \cite{Po}). In the canonical case $D^* = \D^*$, we have the classical Grunsky coefficients.

Consider the set
$$
A_2^2(D) = \{\psi \in A_2(D): \ \psi = \om^2\}
$$
consisting of the integrable holomorphic \ functions \ on $D$ having only zeros of even order and put
$$
\a_D(f) = \sup \ \{|\langle \mu_0, \psi\rangle_D|: \ \psi \in A_2^2,
\ \|\psi\|_{A_2(D)} = 1\}.
$$
The following proposition from \cite{Kr5} completely describes the relation between the Grunsky and Teichm\"{u}ller norms (more special results were obtained in \cite{Kr2}, \cite{Ku2}).

\bigskip\noindent{\bf Proposition 4}. {\em For all $f \in \Sigma^0(D^*)$,
$$
\vk_{D^*}(f) \le k \fc{k + \a_D(f)}{1 + \a_D(f) k}, \quad k = k(f),
$$
and $\vk_{D^*}(f) < k$ unless
 \be\label{11}
\a_D(f) = \|\mu_0\|_\iy,
\end{equation}
where $\mu_0$ is an extremal Beltrami coefficient in the equivalence class $[f]$. The last equality is equivalent to $\vk_{D^*}(f) = k(f)$.

If $\vk_{D^*}(f) = k(f)$ and the class of $[f]$  is a Strebel point, then $\mu_0$ is necessarily of the form}
$$
\mu_0 = \|\mu_0\|_\iy |\psi_0|/\psi_0 \ \ \text{with} \ \ \psi_0 \in A_2^2(D).
$$

\bigskip
Note that geometrically (11) means the equality of the
Carath\'{e}odory and Teichm\"{u}ller distances on the geodesic disk
$\{\phi_\T(t\mu_0 /\|\mu_0\|): t \in \D\}$ in the universal
Teichm\"{u}ller space $\T$.

\bigskip\bigskip
\centerline{\bf 5. PROOF OF THEOREM 4}

\bigskip
We first prove this theorem for $D = \fH$ (and hence for the disk).
In this canonical case, one gets a somewhat stronger result; moreover, the arguments are simpler and illustrate all underlying features.
Now the poles of $r_n$ are real, and $\ld_{\fH}(z) = 1/|z - \ov z| = 1/(2y)$.

The assertion of the theorem follows from the next two lemmas. The first lemma ensures the existence for any $r_n \in \B_2(D)$ of a sequence of points $z_n \in D$ convergent to a boundary point $a_0$ on which the supremum of $\ld_D^{-2}|r_n|$ is attained (such $a_0$
can be distinct from the poles of $r_n$). The second lemma yields that this $a_0$ must be an essential point for $r_n$, and therefore this function represents a non-Strebel point.

Of course, all this is valid to much more general functions from $\B_2(D)$. A special case (the convex hull of fractions $ 1/(z - a)^2$ with real $a$) was considered in \cite{Kr4}.

\bigskip\noindent{\bf Lemma 1}. {\em For any simply connected domain $D \subset \hC$ with quasiconformal boundary $L$ and any rational function $r_n$ with poles of order two on $L$ of the form (3),}
 \be\label{12}
\|r_n\|_{\B_2(D)} = \limsup\limits_{z\to L} \ld_D(z)^{-2} |r_n(z)|. 
\end{equation}

\bigskip\noindent
So, there is a boundary point $z_0$ at which the maximal value of
$\ld_D(z)^{-2} |r_n(z)|$ on $\ov D$ is attained.

\bk\noindent \textbf{Proof}. Consider first the case $D = \D$, and let $r_n \in \B_2(\D)$ satisfy
 \be\label{13}
\limsup\limits_{|z|\to 1} (1 - |z|^2)^2 |r_n(z)| < \sup_{z\in \D} (1 - |z|^2)^2 |r_n(z)|,
\end{equation}
i.e., the polyanalytic function
$$
F(z) = (1 - z \ov z)^2 |r_n(z)|
$$
with $F(0) = r_n(0)$ attains its maximal value on $\ov \D$ at some inner point $z_0 \in \D$. Applying, if needed, the conformal
automorphism
$$
z \mapsto (z - z_0)/(1 - \ov z_0 z)
$$
of $\D$, one reduces the proof to the case $z_0 = 0$.

If $r_n^\prime(0) = a \ne 0$, then $r_n(z) = r_n(0) + a z + \dots$,
and hence, for $z = \rho e^{i \theta}$ and small $\rho > 0$,
$$
\max_\theta |r_n(\rho e^{i \theta})| = |r_n(0)| + |a| \rho +
O(\rho^2).
$$
This yields
$$
\max_\theta |F(\rho e^{i \theta})| = |F(0)| + |a| \rho + O(\rho^2)
> |F(0)|, \quad \rho \to 0,
$$
which contradicts the maximality of $|F(z)|$ at $z = 0$. So, for such rational functions $r_n$,
the inequality (13) can never occur, and 
$$
\limsup\limits_{|z|\to 1} (1 - |z|^2)^2 |r_n(z)| = \sup_{z\in \D} (1 - |z|^2)^2 |r_n(z)| = \|r_n\|_{\B_2(\D)}.  
$$

If $r_n^\prime(0) = 0$, we approximate this function by rational $r_{n,\ve}$ with the same poles $a_k \in \partial \D$, replacing one of the coefficients $c_k$ by $c_k + \ve$ so that $r_{n,\ve}^\prime(0) \ne 0$. 

Since at the point $z_0$, where the function 
$$ 
F_\ve(z) = (1 - z \ov z)^2 |r_{n,\ve}(z)| 
$$ 
attains its maximal value, this value is positive, one can define in a neighborhood of $z_0$ a single valued branch $g_{n,\ve}(z) = \sqrt{r_{n,\ve}(z)}$, and in this neighborhood  
$$
F_\ve(z) = (1 - z \ov z)^2 g_{n,\ve}(z) \ov{g_{n,\ve}(z)}. 
$$
Noting that both partial derivatives  
$\partial_z F_\ve(z), \ \partial_{\ov z} F_\ve(z)$ 
vanish at $z_0$ and 
$$
\partial_z F_\ve(z) = - 2 (1 - z \ov z) \ov z g_{n,\ve}(z) 
\ov{g_{n,\ve}(z)} + (1 - z \ov z)^2 \ov{g_{n,\ve}(z)} 
g_{n,\ve}^\prime(z), 
$$
one obtains   
$$ 
-2 \ov z_0 g_{n,\ve}(z_0) + (1 - z_0 \ov z_0) g_{n,\ve}^\prime (z_0) = 0,   
$$ 
and therefore,  
$$
g_{n,\ve}^\prime(z_0) = - \fc{2 \ov z_0}{1 - z_0 \ov z_0} 
g_{n,\ve}(z) \ne 0. 
$$
This yields, in the same manner as above, that every such $r_{n,\ve}$ satisfies the equality (12). Since this equality remains valid in the limit as $\ve \to 0$, the assertion of Lemma for $D = \D$ is established. 

The case of the generic quasidisk $D$ is reduced to the above one, taking a conformal map $\chi_D$ function of $D$ onto the disk $\D$
with $\chi_D(z_0) = 0$ and applying the above arguments to functions $r_n \circ \chi_D$, which have the same properties as $r_n$. 
This completes the proof of the lemma. 

\bigskip 
Note that $r_n(z) = O(1/z^2)$ as $z \to \iy$, so the quadratic differential $r_n(z) dz^2$ has at the infinite point a pole of the second order. If the boundary of domain $D$ contains $z = \iy$, then the maximal value in (12) can be obtained at this point (and accordingly, $(1 - |\z|^2)^2 |r_n(\z)\chi_D(\z)|$ can attain its maximum at $\z = \chi_D(\iy)$).

\bigskip\noindent{\bf Lemma 2}. {\em Let $D$ be a simply connected
domain on $\hC$ with quasiconformal boundary $L$ and such that $\|S_{g_D}\|_{\B_2} < 1/2$. There exists a constant $c_0 > 0$ such that for any rational function $r_n$ with poles of order two on $L$ of the form (3) and with
$$
\|r_n\|_{\B_2(D)} < c_0 
$$
the boundary points of $D$ at which the maximal value in (12) is attained are substantial for extremal quasiconformal extensions of conformal immersions $f: \ D \to \hC$ generated by the Schwarzian equation $S_f = r_n$ on $D$.}

\bk\noindent\textbf{Proof}. It is sufficient to prove the lemma for domains with boundaries containing $\iy$. We first consider the canonical case $D = \fH$ for which a somewhat stronger result will be obtained.

The equation $S_f(z) = \vp(z)$ defines the conformal immersion
$f_\vp: \fH \to \hC$ determined uniquely by the requirement to preserve the points $0, 1, \iy$.

By the Ahlfors-Weill theorem \cite{AW}, every $\vp \in \B_2(\fH)$ with $\|\vp\| < 1/2$ is the Schwarzian derivative $S_f$ of a
univalent function $f$ in $\fH$, and $f$ has a quasiconformal
extension onto the lower half-plane $\fH^*$ with Betrami coefficient
of the form
 \be\label{14}
\mu_\vp(z) = - 2 y^2 \vp(\ov z), \quad \vp = S_f \ (z = x + i y \in \fH)
\end{equation}
called the {\bf harmonic} Beltrami coefficient (in the spirit of the Kodaira-Spencer deformation theory).

Our aim is to show that {\it for every $r_n$ with real poles $a_j$
of order two and $\|r_n\|_{\B_2(\fH)} < 1/2$ the corresponding harmonic Beltrami coefficient $\mu_{r_n}$ in $\fH$ is extremal in its class, and
 \be\label{15}
\vk(f_{r_n} \circ \sigma) = k(f_{r_n}) = \|r_n\|_{\B_2(\fH^*)},
\end{equation}
where $\sigma$ is the appropriate Moebius map of $\D^*$ onto $\fH^*$.}

It suffices to establish the relations (15) for $r_n$ with
sufficiently small norm.

Pick the point $a_0 \in \R$ at which the equality in (12) is
attained, and two points $x', \ x''$ located on $\R$ in the left to
all poles $a_j$ (so, $\mu_{r_n}(z) = 0$ on $[x^\prime,
x^{\prime\prime}]$). We now establish that
 \be\label{16}
\sup_{\|\psi\|_{A_2(\fH)} = 1} |\langle \mu_{r_n}, \psi
\rangle_{\fH}| = \sup_{\|\psi \|_{A_2^2(\fH)} = 1} |\langle
\mu_{r_n}, \psi \rangle_{\fH}| = H_{a_0}(f),
\end{equation}
which implies the equalities (15) and extremality of $\mu_{r_n}$ in
its class.

Using the conformal map $z = g(\z)$ of the half-strip
$$
\Pi_{+} = \{\z = \xi + i \eta : \ \xi > 0, \ 0 < \eta < 1\}
$$
onto $\fH$ with $g(x^\prime) = 0, \ g(x^{\prime\prime}) = 1, \ g(a_0) = \iy$, we pull-back $\mu_{r_n}/b(a_0)$ (where $b(a_0)$ is the local boundary dilatation at the point $a_0$) to the  Beltrami coefficient
$$
\mu_{*}(\z) := \fc{1}{H_{a_0}} g_{*}(\mu_{r_n})(\z) =
\fc{1}{H_{a_0}}(\mu_{r_n} \circ g)(\z) \
\overline{g^\prime(\z)}/g^\prime(\z)
$$
on $\Pi_{+}$, which satisfies $\lim\limits_{\xi \to \iy}
|\mu_{*}(\xi + i \eta)| = \|\mu_{*} \|_\iy = 1$ and has the limit
function
$$
\mu_{*}(\z_0) = \lim\limits_{\z \to \z_0 \in \partial \Pi_{+}}
\mu_{*}(\z)
$$
with
 \be\label{17}
\mu_{*}(i \eta) = 0.
\end{equation}

We claim that the sequence
$$
\om_m(\z) = \frac{1}{m} e^{- \z/m}, \quad m = 1, 2, \dots \  (\z \in
\Pi_{+}),
$$
is degenerating for $\nu_{*}$. First of all, these $\om_m$ belong to
$A_2^2(\Pi_{+}); \ \om_m(\z) \to 0$ uniformly on $\Pi_{+} \cap
\{|\z| < M\}$ for any $M < \iy$, and $\|\om_m\|_{A_2(\Pi_{+})} = 1$.
Further,
 \be\label{18}
\langle \mu_{*}, \om_m \rangle_{\Pi_{+}} = \fc{1}{m} \iint_{\Pi_{+}}
\mu_{*}(\z) \om_m(\z) d \xi d \eta = \int\limits_0^1 e^{-i \eta/m} d
\eta \ \Bigl( \fc{1}{m} \int\limits_0^\infty \mu_{*}(\xi + i \eta)
e^{- \xi/m} d \xi\Bigr).
\end{equation}
The inner integral can be evaluated using the Laplace transform of $\mu_{*}$ in $\xi$. Integrating by parts and applying (17), one
obtains
$$
\int\limits_0^\iy \fc{\partial \mu_{*}(\xi + i \eta)}{\partial \xi}
e^{- \xi/m} d \xi = \frac{1}{m} \int\limits_0^\iy \mu_{*}(\xi + i
\eta) e^{- \xi/m} d \xi.
$$
On the other hand, Abel's theorem for the Laplace transform yields that the nontangential limit
$$
\lim\limits_{s \to 0} \int\limits_0^\iy \frac{\partial \mu_{*}(\xi +
i \eta)}{\partial \xi} e^{- s\xi} d \xi = \int\limits_0^\iy
\fc{\partial \mu_{*}(\xi + i \eta)}{\partial \xi} d \xi =
\mu_{*}(\iy) - \mu_{*}(i \eta);
$$
hence,
$$
\lim\limits_{m \to \iy} \fc{1}{m} \int\limits_0^\iy \mu_{*}(\xi + i
\eta) e^{- \xi/m} d \xi = \mu_{*}(\iy).
$$
By Lebesgue's theorem on dominated convergence, the iterated
integral in (18) is estimated as follows
 \be\label{19}
\lim\limits_{m \to \iy} |\langle \nu_{*}, \om_m \rangle_{\Pi_{+}}| =
\Big\vert \int\limits_0^1 d \eta \ \lim\limits_{m \to \iy}
 \fc{1}{m} \int\limits_0^\infty \mu_{*}(\xi + i \eta)
e^{- \xi/m} d \xi \Bigr\vert = 1.
\end{equation}
Since by (12), the left-hand side equals to $\|\mu_{*}\|_\iy$ and
all functions $\om_m$ belong to $A_2^2(\Pi_{+})$, this proves our
claim.

Now, applying the inverse conformal map $\z = g^{-1}(z): \Pi_{+} \to
\fH$, one obtains the degenerating sequence
$$
\{\psi_m = (\om_m \circ g^{-1})(g^\prime)^{-2}\} \subset A_2^2(\fH),
$$
for the initial  Beltrami coefficient $\mu_{r_m}$ on $\fH$. By (19),
$$
\lim\limits_{m \to \iy} |\langle \mu_{r_m}, \psi_m \rangle_{\fH}| =
\|\mu_{r_m}\|_\iy = 1,
$$
which implies, together with Lemma 1,  the assertion of Theorem 4
for the half-plane.

\bk Let now $D$ be the generic simply connected domain bounded by
quasicircle $L$ passing through $0, 1, \iy$ and such that the
Schwarzian of the conformal map $g_D$ of $\fH^*$ onto $D$ 
preserving the points $0, 1, \iy$ satisfies
$$\|S_{g_D}\|_{\B_2(\fH^*)} < 1/2. 
$$
Given a rational function $r_n(z)$ of the form (3) with poles on $L$, we consider the univalent solution $w = f_n(z)$ of the equation 
$$
S_w(z) = t r_n(z), \quad z \in D  
$$ 
and its composition with $g_D$, taking $t > 0$ so small that
 \be\label{20}
S_{f_n \circ g_D} = (S_f \circ g_D) (g_D^\prime)^2 + S_{g_D} 
\end{equation}
also satisfies
 \be\label{21}
\|S_{f_n\circ g_D}\|_{\B_2(\fH^*)} < 1/2.
\end{equation}
The  Beltrami coefficients of arbitrary quasiconformal extensions $\wh g_D$ and $\wh f_n$ of $g_D$ and $f_n$, respectively, across the boundaries of their domains to $\hC$ are related by
$$
\mu_{\wh f_n} \circ \wh g_D = \mu_{\wh f_n \circ \wh g_D^{-1}} \circ
\wh g_D = \fc{\mu_{\wh f_n \circ \wh g_D} - \mu_{\wh g_D}}{1 - \ov{\mu_{\wh g_D}} \mu_{\wh f_n \circ \wh g_D}} \ \fc{
\partial_\z \wh g_D}{\ov{\partial_\z \wh g_D}}.
$$
In particular, using their Ahlfors-Weill extensions (14), one gets, in view of (20),
$$
\mu_{\wh f_n} \circ \wh g_D = -2t \eta^2 r_n \circ g_D(\ov \eta) \
\fc{\partial_\z \wh g_D}{\ov{\partial_\z \wh g_D}} + O(t^2)  = - 2t \ld_D^{-2} r_n + O(t^2),
$$
or equivalently, for $z = g_D(\ov \z)$,
$$
\mu_{\wh f}(z) = - 2t \ld_D^{-2}(z) r_n(z) + O(t^2) \quad \text{as}
\ \ t \to 0.
$$
The remainders in the last two equalities are uniformly bounded in
$L_\iy$ norm for all $t$ for which the bound (21) is valid.

We establish now that the harmonic Beltrami coefficient
  \be\label{22}
\nu_{t r_n}(z) = t \ld_D^{-2}(z) r_n(z)
\end{equation}
is {\bf infinitesimally extremal} in its equivalence class.

Indeed, taking again the point $a_0 \in L$ on which the upper limit (12) for chosen $r_n$ 
is attained and mapping the domain $D$ conformally onto the
half-strip $\Pi_{+}$ so that $a_0$ is going to $\iy$, one can repeat
for $f_n$  the above arguments and derive from (16) and (19) that the boundary
dilatation at $a_0$ is equal to $\|r_n\|_{\B_2(D)}$. This yields that the Beltrami coefficient (22) is infinitesimally extremal in $\Belt(D)_1$ and, moreover, its norm is attained on functions from $A_2^2(D)$.

As was mentioned in Section 4.1, such a Beltrami coefficient must be
simultaneously globally extremal in its equivalence class. This implies the assertions of Lemma 2 and of Theorem 4, completing their proofs.

\bigskip
{\small\em{ \leftline{Department of Mathematics, Bar-Ilan
University, 5290002 Ramat-Gan, Israel} 
\leftline{and
Department of Mathematics, University of Virginia, 
Charlottesville, VA 22904-4137, USA}}}

\end{document}